\documentclass[12pt,onecolumn,twoside]{IEEEtran}

\usepackage{graphics} 
\usepackage{epsfig} 
\usepackage{url}
\usepackage{amsmath} 
 \usepackage{amssymb}  
\usepackage{verbatim} 
\usepackage{setspace} 

\newcommand{\bea}{\begin{eqnarray}}
\newcommand{\eea}{\end{eqnarray}}
\newcommand{\beas}{\begin{eqnarray*}}
\newcommand{\eeas}{\end{eqnarray*}}
\newcommand{\be }{\begin{equation}}
\newcommand{\ee }{\end{equation}}

\newcommand{\Q}{\mathcal Q}

\newcommand{\U}{\mathcal U}

\newcommand {\R}{\mathbb R}
\renewcommand{\P}{\mathcal P}

\newcommand{\qed}{\hfill\raisebox{1.2mm}{\fbox{}}}

\newtheorem{Theorem}{Theorem}

\newtheorem{Proposition}{Proposition}

\newtheorem{Example}{Example}
\newtheorem{Remark}{Remark}
  \newtheorem{assumption}{Assumption}
   \newtheorem{conj}{Conjecture}

\newtheorem{Problem}{Problem}

\newcommand{\sname}{}
\newcommand{\slabel}[1]{\debug{\fbox{\tiny \sname #1}}\label{\sname #1}}
\newcommand{\debug}[1]{}              
\newcommand{\FB}{\begin{figure}[t]\centering}
\newcommand{\FE}[2]{\caption{#2
\debug{\fbox{\sname #1}}} \slabel{#1} \end{figure}}
\newcommand{\tB}{\begin{table}[hbtp]\centering}
\newcommand{\tE}[2]{\caption{#2
\debug{\fbox{\sname #1}}}\slabel{#1} \end{table}}

\doublespace
\title{ \LARGE \bf  High-order maximum principles  for the stability analysis of positive  bilinear control~systems}

\author{Gal Hochma$^{1}$  and Michael Margaliot$^{2}$
\thanks{ Research supported in part by the Israel Science Foundation~(ISF). }
\thanks{$^{1}$GH is with the School of Electrical Engineering-Systems, Tel
Aviv University, Israel 69978. Email: \tt{gal.hoc@gmail.com}}
\thanks{$^{2}$MM (corresponding author) is with the School of Electrical Engineering-Systems
and with the Sagol School of Neuroscience, Tel
Aviv University, Israel 69978. Email: \tt{michaelm@eng.tau.ac.il}}
\thanks{An abridged version of this paper was presented at the
\emph{$52$nd IEEE Conference on Decision and Control}~\cite{sec_order_cdc}. }
}

\begin{document}

\maketitle
\thispagestyle{empty}
\pagestyle{empty}

\begin{abstract}
 We consider a continuous-time positive bilinear control system~(PBCS), i.e.
  a bilinear control system with
Metzler matrices.
The positive orthant  is an invariant set of such a system, and
 the corresponding transition matrix~$C(t)$  is
entrywise nonnegative for all time~$t\geq 0$.
Motivated by the stability analysis
of positive linear switched systems~(PLSSs) under arbitrary switching laws, we fix
a final time~$T>0$
and define a control as
optimal if  it maximizes the spectral radius of~$C(T)$.
A recent paper~\cite{lior_SIAM} developed
  a first-order necessary condition for optimality
  in the form of  a maximum
principle~(MP). In this paper, we derive   higher-order  necessary conditions for optimality for
 both singular and bang-bang controls.
Our approach is based on combining results on the second-order derivative of the spectral radius of a nonnegative
matrix with  the generalized Legendre-Clebsch condition and
the Agrachev-Gamkrelidze second-order optimality condition.
\end{abstract}



\section{Introduction}\label{sec:intro}
Consider the continuous-time  linear switched  system
\begin{align}\label{eq:a0a1}
\dot{x}(t)& =A_{\sigma (t)}x(t),   \nonumber \\
x(0)& =x_{0},
\end{align}%
where~$x:\mathbb{R}_{+}\rightarrow \mathbb{R}^{n}$ is the state vector, and~$%
\sigma :\mathbb{R}_{+}\rightarrow \{0,1\}$ is a piecewise constant function
referred to as the \emph{switching signal}. This models a system that can
switch between the two linear subsystems
\begin{equation*}
\dot{x}=A_{0}x\;\text{ and }\;\dot{x}=A_{1}x.
\end{equation*}

Recall that~\eqref{eq:a0a1} is said to be  \emph{globally uniformly asymptotically stable}~(GUAS) if there exists a class~$\mathcal{K}\mathcal{L}$ function\footnote{A continuous function $\alpha :[0,\infty)\rightarrow [ 0,\infty )$
belongs to the class $\mathcal{K}$ if it is strictly increasing and
$\alpha (0)=0$. A continuous function $\beta :[0,\infty)\times [ 0,\infty
)\rightarrow [ 0,\infty )$ belongs to the class $\mathcal{K}\mathcal{L}
$ if for each fixed $s$, $\beta (\cdot ,s)$ belongs to $\mathcal{K}$, and
for each fixed $r>0
$, the mapping $\beta (r,\cdot )$ is decreasing and $\beta
(r,s)\rightarrow 0$ as $s\rightarrow \infty $.}~$\beta $ such that for any
initial condition~$x_0 \in \mathbb{R}^{n}$ and \emph{any} switching law~$%
\sigma $, the corresponding solution of~\eqref{eq:a0a1} satisfies
\begin{equation*}
|x(t)|\leq \beta (|x_0|,t),\text{ for all }t\geq 0.
\end{equation*}%
This implies in particular that
\begin{equation}
\lim_{t\rightarrow \infty }x(t)=0,\quad \text{ for all }\sigma
\text{ and all }x_0 \in \mathbb{R}^{n},  \label{eq:ctz}
\end{equation}%
and or linear switched systems,~\eqref{eq:ctz} is in fact equivalent to GUAS (see,
e.g.,~\cite{angeli-sontag-positive-2009}).
Switched systems and, in particular, their stability analysis
are attracting considerable interest  in the last two decades;
see e.g. the survey papers~\cite{shorten,branicky_98,decarloetal,libsur99} and the monographs~\cite{liberzon_book,morse-book,Jungers,johansson,sun_ge_book, sun_ge_2}.

It is well-known   that a
necessary (but not sufficient) condition for GUAS of~\eqref{eq:a0a1} is the following.
\begin{assumption} \label{a:nec}
The matrix~$kA_{0}+(1-k)A_{1}$ is
Hurwitz for all~$k\in \lbrack 0,1]$.
\end{assumption}

Recall that a linear system
\be\label{eq:posys}
\dot{x}=Ax,
\ee
 with~$A\in \R^{n\times n}$,
is called \emph{positive} if the positive orthant
\begin{equation*}
\mathbb{R}_{+}^{n} : =\{x\in \mathbb{R}^{n}\mid x_{i}\geq 0,\;i=1,\dots, n\}
\end{equation*}%
is an invariant set of the dynamics, i.e.,~$x(0)\in \mathbb{R}_{+}^{n}$
implies that~$x(t)\in \mathbb{R}_{+}^{n}$ for all~$t\geq 0$.

Positive   systems play an important role in
systems and control theory because in many physical systems the
state-variables represent quantities that can never attain negative values
(e.g. population sizes,
 probabilities,      concentrations,   buffer loads)~\cite{farina2000,berman87,posi_sys89}.
A necessary and
sufficient condition for~\eqref{eq:posys} to be positive
  is that $A$ is a \emph{Metzler matrix}, that
is, $a_{ij} \geq 0 $ for all~$i \not = j$. If~$A$ is Metzler
then~$\exp (At)$ is (entrywise)   nonnegative
for all~$t\geq 0$. By the
Perron--Frobenius theory, the spectral radius of~$\exp (At)$ (i.e., the eigenvalue   with maximal
absolute value) is real and nonnegative, and since~$\exp(At)$ is non-singular, it is in fact positive.

If both~$A_{0}$ and~$A_{1}$ are Metzler and~$x(0)\in \mathbb{R}_{+}^{n}$
then~\eqref{eq:a0a1} is called  a \emph{positive linear switched system}%
~(PLSS).
Mason and Shorten~\cite{mason-shorten03}, and independently David Angeli,
 posed the following.
\begin{conj}\label{conj:shaorten}
\label{conj} If~\eqref{eq:a0a1} is a PLSS, then Assumption~\ref{a:nec}
provides a \emph{sufficient} condition for GUAS.
\end{conj}

Had this conjecture been true, it would have implied that
 determining GUAS for  a PLSS is relatively simple.
 (See~\cite{Gurvits_Olshevsky} for analysis of the computational
 complexity of  determining whether any matrix in a convex set of matrices is Hurwitz.)
 Gurvits, Shorten, and Mason~\cite{gurvits-shorten-mason07} proved that Conjecture~\ref{conj:shaorten}
  is in general false (see also~\cite{gurvits-cdc03}), but 
 that it does hold when~$n=2$ (even when the number of
subsystems is arbitrary). Their proof in the planar case is based on showing
that the PLSS admits a common quadratic Lyapunov function~(CQLF). (For more on
the analysis of switched systems using CQLFs, see~\cite%
{branicky_98,shorten,copos,forna,ron_margaliot}.) Margaliot and Branicky~\cite{mar-bra-full}
  derived a
reachability--with--nice--controls--type result for planar bilinear control
systems, and showed that the proof of Conjecture~\ref{conj} when~$n=2$
follows as a special case. Fainshil, Margaliot, and Chigansky~\cite{lior} showed 
 that Conjecture~\ref{conj} is false already for the case~$n=3$. 
 In general, it seems that   as far as the GUAS problem is concerned,
analyzing PLSSs is  not simpler than analyzing linear switched
systems.

There is a rich literature  on \emph{sufficient} conditions for GUAS, see, e.g., \cite{branicky_98,decarloetal,liberzon_book,libsur99,sun_ge_book}.
A more challenging problem is to determine a \emph{necessary and sufficient} condition for GUAS.
What makes this
 problem    difficult is that  the set of all possible switching laws is
huge, so exhaustively checking the solution for each switching law is
impossible.

 A natural idea  is to try and characterize a ``most destabilizing''
  switching law~$\sigma ^{\ast }$ of the switched system, and
then analyze the behavior of the corresponding trajectory~$x^{\ast }$. If~$%
x^{\ast }$ converges to the origin, then so does any trajectory of the
switched system and this establishes GUAS.
This idea was    pioneered by E. S.
Pyatntisky~\cite{pyat70,pyat71},
who   studied the celebrated
\emph{absolute stability problem}~(ASP).
 This
variational approach was further developed by several scholars including
 N. E. Barabanov
and L. B. Rapoport, and proved to be highly successful; see the
survey papers~\cite{mar-simple,bar-cdc,rap}, the related work in~\cite{Boscain2008Nondiagonalizable_case,Boscain2009stability_conditions},
and
 the recent extensions to the stability analysis of
 {discrete--time} linear switched systems in~\cite{monovich1,tal2}.

A  first attempt to extend the variational approach to the stability
  analysis of PLSSs was taken
  in~\cite{mar-bra-posi} using
 the classical Pontryagin maximum principle~(PMP).
Recently,  Fainshil and Margaliot~\cite{lior_SIAM} developed  an alternative
 approach
that combines the Perron-Frobenius theory of nonnegative matrices with
the standard needle variation used in the PMP.

The goal of this paper is to derive stronger, higher-order necessary conditions for optimality.
We thus begin by reviewing the first-order MP in~\cite{lior_SIAM}.

\subsection{Stability analysis of PLSSs: a Variational  Approach}
The  variational  approach to the stability analysis of a linear switched system includes several steps.
  The first step
 is relaxing~\eqref{eq:a0a1} to the \emph{bilinear control system}~(BCS)
 \begin{align}\label{eq:pscon}
\dot{x}& =(A+uB)x,\quad u\in \mathcal{U},   \\
x(0)& =x_{0}\nonumber,
\end{align}
where~$A: =(A_0+A_1)/2$~,~$B:=(  A_1-A_0)  /2$, and~$\U$ is the set of measurable controls taking values in~$[-1,1]$. Note that for~$u(t) \equiv -1$ [$u(t) \equiv 1$],
Eq.~\eqref{eq:pscon} yields~$\dot{x}=A_0 x$ [$\dot{x}=A_1x]$, i.e.,
trajectories of the~BCS corresponding to
piecewise constant  bang-bang controls are
also trajectories of the original switched system.

The BCS~\eqref{eq:pscon}  is said to be \emph{globally asymptotically stable}~(GAS) if~$\lim_{t \to \infty} x(t)=0$ for all~$x_0\in\R^n$ and \emph{all}~$u \in \mathcal{U}$.
Since every trajectory of the switched system~\eqref{eq:a0a1} is also a trajectory of~\eqref{eq:pscon},
GAS of~\eqref{eq:pscon} implies GUAS of the linear switched system.
It is not difficult to show that the converse implication also holds, so the BCS is GAS
if and only if the linear switched system is GUAS.
Thus, the GUAS problem for the switched linear system~\eqref{eq:a0a1} is equivalent to the GAS
problem for the BCS~\eqref{eq:pscon}.

From here  on we assume that the switched system is positive, 
i.e.~$A+kB$ is Metzler for all~$k\in[-1,1]$. For the BCS, 
this implies that if~$
x_{0}\in \mathbb{R}_{+}^{n}$, then~$x(t)\in \mathbb{R}_{+}^{n}$
for all $u\in \mathcal{U}$\ and all~$t\geq 0$.
Thus~\eqref{eq:pscon} becomes a \emph{positive bilinear control system}~(PBCS).

For~$0\leq a \leq b \leq     T$, and~$u \in \U$,
let~$C(b,a,u)$ denote the solution  at time~$t=b$
of the matrix differential equation
\begin{align} \label{eq:trans}
\frac{d}{dt}C(t,a,u)& =(A+Bu(t))C(t,a,u ),  \notag
\\
C(a ,a,u )& =I.
\end{align}
It is straightforward to verify that  the solution  of~\eqref{eq:pscon}
satisfies~$
x(b)=C(b,a,u)x(a)$
for all~$u \in \U$
and     all~$0\leq a \leq b \leq     T$. In other words,~$C(b,a,u)$ is the \emph{transition matrix}
from time~$a$ to time~$b$ of~\eqref{eq:pscon} corresponding to the
control~$u$. To simplify the notation, we will sometimes omit the dependence
  on~$u$ and just write~$C( b,a)$.

When the initial time is~$a=0$ we write~\eqref{eq:trans} as
\begin{align} \label{eq:transzero}
\dot{C}(t)& =(A+Bu(t))C(t),  \notag
\\
C(0 )& =I.
\end{align}
For a PBCS,~$C(t,u)$ is a non-negative matrix for all~$t\geq 0$ and all~$u\in \mathcal{U}$. Since it
 is also
 non-singular, the spectral radius~$\rho
(C(t,u))$ is a real and positive eigenvalue of~$C(t,u)$, called the Perron root.
If this eigenvalue is simple then the corresponding eigenvector~$v \in \R^n_+$, called the Perron eigenvector,
is unique (up to multiplication by a scalar).
The next step in the variational approach is to relate~$\rho
(C(t,u))$ to GAS of the PBCS.

Define the \emph{generalized spectral radius}
 of the PBCS~\eqref{eq:pscon} by
\[
                    \rho(A,B):=\limsup_{t \to \infty} \rho_t(A,B),
\]
where
\be \label{eq:rhot}
\rho_t(A,B):=   \max_{u \in \U} (\rho(C(t,u)))^{1/t}     .
\ee
Note that the maximum here is well-defined, as the reachable set of~\eqref{eq:transzero}
corresponding to~$\mathcal{U}$ is compact~\cite{filippov-paper}.
In fact, this is why we    consider a bilinear control system with controls in~$ \U$
rather than the original linear switched system with piecewise constant switching laws.

  The next result relates the GAS of the PBCS to~$\rho(A,B)$.
\begin{Theorem} \label{thm:starho}
               The PBCS~\eqref{eq:pscon} is GAS if and only if \[\rho(A,B)<1.\]
\end{Theorem}

Thm.~\ref{thm:starho} already appeared in~\cite{lior_SIAM}, but without a proof.
For the sake of completeness we include its proof in the Appendix.

\begin{Remark}
It follows from~\eqref{eq:rhot}  and Thm.~\ref{thm:starho}
that if~$\rho(C(T ,u)) \geq 1$ for some~$T>0$ and~$u \in \U$,
then the~PBCS is \emph{not} GAS.
Indeed,
for any integer~$k>0$, define~$\bar{u}:[0,kT] \to [-1,1]$ via the periodic extension of~$u$,
and let~$\bar{C}(t)$ denote the corresponding solution of~\eqref{eq:transzero} at time~$t$. Then
\[
\rho( \bar{C}( kT ) )=(\rho(C(T)))^k  ,
\]
so~\eqref{eq:rhot} yields
\begin{align*}
\rho_{kT}(A,B)&\geq  (\rho( \bar{C}(kT)))^{1/(kT)}\\&=(\rho(C(T )))^{1/T}\\& \geq 1,
\end{align*}
and this implies that~$\rho(A,B) \geq 1$.~\qed
\end{Remark}

 Thm.~\ref{thm:starho} motivates
the following optimal control problem.
\begin{Problem}
\label{prob:rho} Consider the PBCS~\eqref{eq:transzero}. Fix
an arbitrary final time~$T>0$. Find a control~$u^{\ast }\in \mathcal{U}$ that \emph{maximizes}~$\rho (C(T ,u)) $.
\end{Problem}

The main result in~\cite{lior_SIAM} is a first-order necessary  condition  for optimality. Let~$A'$
  denote the transpose of the matrix~$A$.

\begin{Theorem}\label{thm:mprho}\cite{lior_SIAM}
Consider the PBCS~\eqref{eq:transzero}. Suppose that~$u^*\in \mathcal{U}$ is an optimal control for
Problem~\ref{prob:rho}. Let~$C^*(t)  $ denote the
corresponding solution of~\eqref{eq:transzero} at time~$t$, and let~$\rho ^{\ast }:=\rho
(C^*(T))$.  Suppose that~$\rho ^{\ast }$ is  a \emph{simple} eigenvalue of~$C^*(T)$.
Let~$v^{\ast }\in \mathbb{R}_{+}^{n}$ [$w^{\ast }\in \mathbb{R}_{+}^{n}$]
 be an  eigenvector of~$C^*(T)$ [$(C^*(T))'$] corresponding
to~$\rho^*$, normalized such that~$(w^*)'v^*=1$.
   Let~$q:[0,T]\to \mathbb{R}_{+}^{n}$ be the solution of
\begin{align}\label{Qeq}
\dot{q}& =-(A+Bu^{\ast })'q,   \\
q(T)& =w^{\ast },  \notag
\end{align}%
and let~$p:[0,T]\to \mathbb{R}_{+}^{n}$ be the solution of
\begin{align}
\dot{p}& =(A+Bu^{\ast })p,  \label{Peq} \\
p(0)& =v^{\ast }.  \notag
\end{align}%
Define the \emph{switching function}~$m:[0,T]\to \mathbb{R}$ by
\begin{equation}
m(t):=q'(t)Bp(t).  \label{SwitchingFunction}
\end{equation}%
Then for almost all $t\in \left[ 0,T\right]$,
\begin{equation}\label{eq:UStar}
u^{\ast }(t)=%
\begin{cases}
1, & m(t)>0, \\
-1, & m(t)<0.%
\end{cases}
\end{equation}
\end{Theorem}

This MP has some special properties.
\begin{Remark}
First, note that~\eqref{Qeq} implies that
\[
                  q'(0)=  q'(t)C^*(t), \quad \text{ for all } t \in[0,T].
\]
In particular, substituting~$t=T$ yields
\begin{align*}
                    q'(0)&=q'(T) C^*(T) \\
                          &=(w^*)'C^*(T) \\
                          &=\rho^* (w^*)' ,
\end{align*}
as~$w^*$ is an eigenvector of~$(C^*(T))'$ corresponding to the eigenvalue~$\rho^*$.
Since  scaling~$q$ by a positive constant has no effect on the sign of~$m$, this means
that
  the final condition~$q(T)=w^*$ in~\eqref{Qeq} can be replaced by
  the initial  condition~$q(0)= w^* $. This leads to
   an MP in the form of    a \emph{one-point} boundary value problem
  (with the unknown~$v^*,w^*$ as the initial  conditions at time~$0$).~\qed
\end{Remark}

\begin{Remark} Note that
\begin{align}\label{eq:misperiodic}
                    m(T)&=q'(T) B p(T) \nonumber \\
                        &=(w^*)' B C^*(T) p(0)\nonumber  \\
                        &=(w^*)' B \rho^* v^* \nonumber \\
                        &=q'(0) B  p(0)\nonumber  \\
                        &=m(0).
\end{align}
Thus, the switching function is ``periodic'' in the sense that~$m(T)=m(0)$.~\qed
\end{Remark}

One difficulty in applying Theorem~\ref{thm:mprho}
is that both~$v^*$ and~$w^*$ are unknown.
There are cases where this difficulty may be alleviated somewhat
berceuse~$w^*$ can be expressed in terms of~$v^*$.
 The next example demonstrates this.
 \begin{Example}\label{exa:w_fiunc_v}
 Consider an optimal bang-bang
 control in the form
 \[
            u^*(t)=\begin{cases}
            1,& t \in (0,t_1),\\
            -1,& t \in (t_1,T),\end{cases}
 \]
 where~$0<t_1<T$.
 The corresponding transition matrix is
 \[
            C^*(T)=\exp((A-B)\tau_2)\exp((A+B)\tau_1),
 \]
 where~$\tau_1:=t_1-0$ and~$\tau_2:=T-t_1$.
 Thus,~$v^*$ and~$w^*$ satisfy
 \begin{align*}
                                   \exp((A-B)\tau_2)\exp((A+B)\tau_1)v^*&= \rho^* v^*,
 \end{align*}
  and
   \begin{align}\label{eq:wmulle}
                                   \exp((A+B)'\tau_1)\exp((A-B)'\tau_2)w^*&= \rho^* w^*.
    \end{align}
 Suppose that~$A$ and~$B$ are symmetric matrices. Then~\eqref{eq:wmulle} becomes 
 \[
 \exp((A+B) \tau_1)\exp((A-B) \tau_2)w^* = \rho^* w^*,
 \]
 and multiplying this on the left by~$\exp((A-B) \tau_2)$ yields
 \[
C^*(T) \exp((A-B) \tau_2)w^* = \rho^* \exp((A-B) \tau_2) w^*.
 \]
 Since the Perron eigenvector of~$C^*(T)$ is unique (up to   multiplication by a constant) this means that
 \[
 \exp((A-B) \tau_2)w^* = r v^*,
 \]
 for some~$r >0$.~\qed
 \end{Example}

The MP in Theorem~\ref{thm:mprho}   is a necessary, but not sufficient, condition for optimality and it is possible of course
that
a
 control satisfying this MP is not an optimal control.
 The next example demonstrates  this.
\begin{Example}\label{exa:mp2}
Consider a  PBCS satisfying
the following properties:
\begin{itemize}
\item  The matrix $A $ is symmetric. Its
   maximal eigenvalue~$\mu$  is simple with corresponding eigenvector~$z$, and
\be\label{eq:vpbv}
z ' B z=0.
\ee
\item  The matrices~$A-B$ and~$A+B$ are Metzler;
\item  $\rho(A+B)>\rho(A)=\mu$.
 \end{itemize}
(A specific example is~$n=2$,
$A=\begin{bmatrix}
 2.2 &  1.6 \\ 1.6  &  -0.2
\end{bmatrix}$
and
$B=\begin{bmatrix}
 -1.1   &  0.2 \\
  0.95  & 2.1
\end{bmatrix}
$.
Indeed, here~$\mu=3$, $z = \begin{bmatrix}2&1\end{bmatrix} '$ and it is straightforward to verify that all
the properties above hold.)

Consider the possibility that the singular control~$u (t)  \equiv  0$ is optimal.
Then
\[
\rho(u):=\rho(\exp(AT))=\exp(\mu T).
\]
 Since~$A$ is symmetric,
the corresponding right and left eigenvector is~$z$, so in the MP $p(0)=q(T)=z$.
Thus,~\eqref{Peq} and~\eqref{Qeq} yield
\begin{align*}
p(t)&= \exp(At)z\\
    & =\exp(\mu t) z,
\end{align*}
and
\begin{align*}
q(t)&= \exp(A'(T-t))z\\
    & = \exp(\mu(T-t)) z.
\end{align*}
Substituting this in~\eqref{SwitchingFunction} yields
\begin{align*}
m(t)&=\exp(\mu T)z'   B z \equiv 0.
  \end{align*}
Thus,~$u (t)\equiv 0$  (vacuously) satisfies Thm.~\ref{thm:mprho}.
However, since~$\rho(A+B)>\rho(A)$ the control~$\tilde u(t)\equiv 1$ yields
\[
        \rho(\tilde u):=\exp((A+B)T)>\rho(u),
\]
so clearly~$u(t)\equiv 0$ is \emph{not} an optimal control.~\qed
\end{Example}

The reason that~$u(t)\equiv 0$ in Example~\ref{exa:mp2}
 cannot be ruled out  is that   Thm.~\ref{thm:mprho} is a \emph{first-order} MP.
More specifically, its derivation is based the following idea. Suppose that~$u$
is a candidate for an optimal control.
Introduce a new control~$\tilde{u}:[0,T] \to [-1,1]$ by adding a needle variation to~$u $, i.e.
\[
        \tilde{u}(t):=\begin{cases} a , &t \in [\tau,\tau+\epsilon),\\
                            u (t),& \text{otherwise},\end{cases}
\]
where~$a \in[-1,1]$,~$\tau \in [0,T)$ is a Lebesgue point of~$u $,
and~$\epsilon>0$ is sufficiently small, and analyze
the difference~$C(T,\tilde{u})-C(T,u^*)$ to first-order in~$\epsilon$.
For~$u (t) \equiv 0$,
\[
            C(T, \tilde{u})=\exp(A  (T-\tau-\epsilon ) ) \exp( (A+a B)  \epsilon) \exp(A  \tau),
\]
so
\begin{align*}
                     \left.    \frac{d}{d \epsilon}    C(T,\tilde{u})\right |_{\epsilon=0} & =     a \exp(A  (T-\tau ) )   B
                                 \exp(A  \tau)  .
\end{align*}
Combining this with
  known results on the derivative of a simple eigenvalue of a matrix (see, e.g.~\cite[Chapter~6]{mat_ana_sec_ed}) yields
 \begin{align}\label{eq:oeps}
           \rho(C&(T, \tilde{u}))    = \rho(C(T,u ))   +\epsilon  a w'  \exp(A  (T-\tau  ) )   B
                                \exp(A  \tau)v    +o(\epsilon)  .
\end{align}
If~$a w'  \exp(A  (T-\tau  ) )   B \exp(A  \tau)v  >0$ then~$\rho(C(T,\tilde{u}))    > \rho(C(T,u ))$
for all sufficiently small~$\epsilon>0$ and thus~$u $ is not optimal. However,
in Example~\ref{exa:mp2}   the term multiplying~$\epsilon$ in~\eqref{eq:oeps} is zero
for all~$ a$, $\tau$, and $T$,   so a first-order analysis cannot rule out the
possibility that~$u $ is optimal.

Summarizing, Example~\ref{exa:mp2} suggests that there is a need for a higher-order
 MP, i.e.,
 an MP that takes into account higher-order terms in the Taylor expansion
of~$\rho(C(T,\tilde{u}))   - \rho(C(T,u ))$ with respect to~$\epsilon$,
and can thus be used
to rule out  the optimality of a larger set of controls.

In the next section, we
apply the  generalized Legendre-Clebsch condition to derive
a high-order necessary condition  for a singular control to be optimal.
We also
combine known results on the  second-order  derivative of the Perron root~\cite{Deutsch19841}
and the
Agrachev-Gamkrelidze second-order variation
 for bang-bang controls (see, e.g.,~\cite{agrachev-sigalotti}) to derive a \emph{second-order} MP for
  bang-bang controls.
The proofs of these results are given in Section~\ref{sec:proofs}.


\section{Main results}
Our first result is a  high-order
 necessary condition
  for singular optimal controls for Problem~\ref{prob:rho}. Without loss of generality (see~\cite{hermes78}),
 we assume that the singular control is~$u^*(t) \equiv 0$.
Let~$[P,Q]:=QP-PQ$   denote the  Lie-bracket
of   the matrices~$P,Q \in \R^{n \times n}$.
\subsection{High-order MP for singular controls}

\begin{Theorem}\label{thm:secsing}
Consider the PBCS~\eqref{eq:pscon}. Suppose that the conditions
  of Thm.~\ref{thm:mprho} hold, and
    that~$u^*(t) \equiv 0$ is an optimal control.
 Then
\begin{equation}\label{eq:abb}
(w^*)'[B,[B,A]]  v^* \leq 0.
\end{equation}
\end{Theorem}

\begin{Example}\label{exa:sing_sec}
Consider  the specific PBCS with~$n=2$ given  in Example~\ref{exa:mp2}.
In this case,
\[
[B,[B,A]]  =\begin{bmatrix} 6.8& 18.4\\ 21.4& -6.8\end{bmatrix},
\]
and~$v^*=w^* = \begin{bmatrix}2&1\end{bmatrix} '/\sqrt{5}$,
so
\[
      (w^*)'[B,[B,A]] v^* = 20.
\]
It follows from~\eqref{eq:abb} that~$u^*(t)\equiv 0$ is \emph{not} an optimal control.
Note that we were not able to derive this conclusion using the first-order MP
in Thm.~\ref{thm:mprho}.~\qed
\end{Example}

\subsection{Second-order MP for bang-bang controls}
In this section, we derive an Agrachev-Gamkrelidze-type second-order MP for optimal bang-bang controls for Problem~\ref{prob:rho}.
Note that for an optimal bang-bang~$u^*$ we have
\[
            C^*(T)=\exp((A+B)\tau_k)  \dots\exp((A+B)\tau_2) \exp((A-B)\tau_1) ,
\]
with~$\tau_i\geq 0$ and~$\sum_{i=1}^k \tau_i=T$. Any cyclic shift of~$C^*(T)$, e.g.,
\[
         \exp((A-B)\tau_1)  \exp((A+B)\tau_k)  \dots\exp((A+B)\tau_2)
\]
 also corresponds to an optimal control
(as a product of matrices and its cyclic shift have the same spectral radius).
This means that we can always assume that~$t_0:=0$ is a switching point of~$u^*$, and then~\eqref{eq:misperiodic}
implies that~$T$ is also a switching point of~$u^*$.

Let~$ \P^k$ denote the set of all vectors~$\alpha=\begin{bmatrix} \alpha_0 & \alpha_1&\dots& \alpha_k\end{bmatrix}'\in\ \R^{k+1}$
satisfying
\be \label{eq:defpk}
 \alpha_1 + \dots + \alpha_k = 0  .
 \ee
   We can now state the main result in this section.
\begin{Theorem}\label{thm:mainhere}
Suppose that~$u^*$ is an optimal control for Problem~\ref{prob:rho},
that the conditions of Thm.~\ref{thm:mprho} hold, and that   the   switching function~\eqref{SwitchingFunction} admits
a finite number of zeros at~$  t_0 < t_1 < \dots < t_k  $, with~$t_0=0$,~$t_k= T$,
so that~$u^*(t)=r$ for~$t \in (0,t_1)$, $u^*(t)=-r$ for~$t \in ( t_1,t_2)$,
$u^*(t)=r$ for~$t \in ( t_2,t_3)$, and so on,
with~$r \in \{-1,1\}$. Denote~$P:=A+rB$,~$Q:=A-rB$,
and~$\tau_i:=t_i-t_{i-1}$.
Define matrices~$H_i \in \R^{n\times n}$, $i=1,\dots,k $, by
\begin{align}\label{eq:defhis}
                        H_1&:=P,\\
                        H_2&:=Q,\nonumber\\
                        H_3&:=\exp(-\tau_2 Q) P \exp(\tau_2 Q),\nonumber\\
                        H_4&:=\exp(-\tau_2 Q) \exp(-\tau_3 P)Q \exp(\tau_3 P) \exp(\tau_2 Q),\nonumber\\
                        H_5&:=\exp(-\tau_2 Q) \exp(-\tau_3 P) \exp(-\tau_4 Q) P \exp(\tau_4 Q)  \exp(\tau_3 P) \exp(\tau_2 Q)\nonumber,\\
                           &\vdots\nonumber
\end{align}
Then
\be \label{eq:firsth1}
                           q'(t_1) \sum_{i=1}^{k } \alpha_i  H_i p(t_1)=0, \quad \text{for all } \alpha \in \P^{k }.
\ee
Furthermore,
\[
            r_k(\alpha):=q'(t_1) \sum_{1\leq i < j \leq k } \alpha_i \alpha_j [H_i,H_j] p(t_1)
\]
satisfies
\be \label{eq:secoopt}
                r_k(\alpha) \leq 0, \quad \text{for all } \alpha \in \Q^{k },
\ee
where
\be\label{eq:defpkpk}
                \Q^{k } :=\{ \alpha \in \P^{k }  : \sum_{i=1}^{k } \alpha_i H_i p(t_1) =0  \}.
\ee
\end{Theorem}

We refer to the control~$u^*$ defined above as a   control with~$k$ bang arcs. As will be shown in the proof, condition~\eqref{eq:firsth1} is a first-order condition (that can also be derived using  the first-order MP).
Condition~\eqref{eq:secoopt} however is a second-order condition, and it is meaningful for values~$\alpha$
that make a certain  first-order variation vanish, i.e. that belong to~$\Q^k$.

Note that the conditions in Thm.~\ref{thm:mainhere} are given in terms of~$p(t_1)$ and~$q(t_1)$. It is possible of course
to state them in terms of~$p(t_0)=v^*$ and~$q(t_0)= \rho^* w^*$, but this leads to
slightly more cumbersome expressions.

\begin{figure}[t]
\begin{center}
\mbox{\psfig{figure=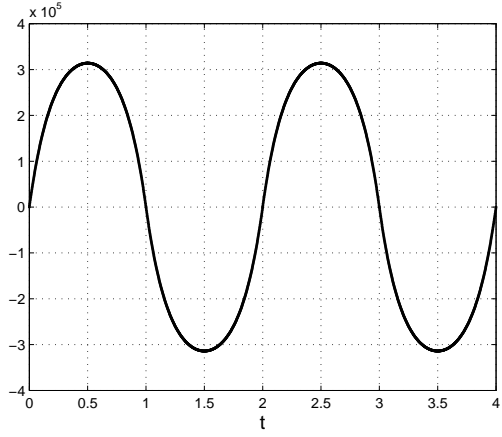,width=8cm,height=6cm}}
\end{center}
\caption{Switching function $m(t)$ in Example~\ref{exa:2dsimpl}.  }
\label{fig:m4s}
\end{figure}

The next example demonstrates the calculations for a  control with two bang arcs.
\begin{Example}\label{exa:1switch}
 Consider an optimal control in the form
 \[
            u^*(t)=\begin{cases}
            1,& t \in (0,t_1),\\
            -1,& t \in (t_1,T),\end{cases}
 \]
 where~$0<t_1<T$.
 In this case,~\eqref{eq:firsth1} becomes
\[
     q'(t_1) (\alpha_1 (A+B) + \alpha_2 (A-B))   p(t_1)=0, \quad \text{for all } \alpha \in \P^{2 },
\]
and the definition of~$\P^{2 }$ yields
\[
    \alpha_1 q'(t_1) (  (A+B) -   (A-B))   p(t_1)=0, \quad \text{for all } \alpha_1 \in \R.
\]
Of course, this is just the conclusion that we can get from the first-order MP, as at the switching point~$t_1$ we must have
\[
        0=m(t_1)=q'(t_1) B  p(t_1).
\]
The second-order term is
\begin{align*}
            r_2(\alpha) &=\alpha_1 \alpha_2 q'(t_1)     [H_1,H_2] p(t_1)\\
                        &=- \alpha_1^2 q'(t_1)    [A+B,A-B] p(t_1)\\
                        &=  2\alpha_1^2 q'(t_1)    [ A,B ] p(t_1),
\end{align*}
so~\eqref{eq:secoopt} becomes
\be\label{eq:r2eq}
                r_2(\alpha) \leq 0, \quad \text{for all } \alpha \in \Q^{2 },
\ee
where
\[
                \Q^{2 }   =\{ \alpha_1 \in \R  : \alpha_1 B  p(t_1) =0  \}.
\]
Again, this provides information that can also be derived from the first-order MP, as the fact that~$m(t_1^-)>0$ and~$m(t_1^+)<0$ implies that
\[
        \dot m(t_1)  \leq 0.
\]
and differentiating~\eqref{SwitchingFunction} yields
\[
\dot m(t_1) =  q'(t_1)[A,B] p(t_1) .
\]
Thus,~$q'(t_1)[A,B] p(t_1) \leq 0$, so~\eqref{eq:r2eq} actually holds for all~$\alpha \in \P^2$.~\qed
\end{Example}
However, for a control with more than two bang arcs the second-order condition does
provide  new information.
The next simple example demonstrates this.
\begin{Example}\label{exa:2dsimpl}
Consider the PBCS~\eqref{eq:transzero} with
\begin{align*}
A=\begin{bmatrix} -5/2 & 3/2 \\ 3 & -5/2\end{bmatrix},
\quad
B=\begin{bmatrix} 3/2 & -1/2 \\ 1 & -3/2\end{bmatrix}.
\end{align*}
Note that~$A+kB$ is Metzler for all~$k\in[-1,1]$. Consider the control
 \be\label{eq:opcande}
            u (t)=\begin{cases}
            1,& t \in (0,t_1) \cup (t_2,t_3), \\
            -1,& t \in (t_1,t_2) \cup (t_3,T),
            \end{cases}
 \ee
 with~$t_1=1$, $t_2=2$, $t_3=3$, and~$T=4$.
 The corresponding transition
 matrix is
 \[
 C(T)=\exp(A-B) \exp(A+B)\exp(A-B) \exp(A+B).
 \]
 Let~$s:=\left( 9+32\exp( 5 )+9\exp( 10 ) \right)^{1/2}$.
 The spectral radius  of~$C(T)$  is
 \[
            \rho=\left (  \frac{ 9+7\exp(5)+9\exp(10)+3s(\exp(5)-1)   }{  25\exp(10)  }   \right)^2 ,
 \]
and it is  a simple eigenvalue.
The Perron right and left eigenvectors
 of~$C(T)$  are
 \[
 v=\begin{bmatrix}
                                         \exp(5)-1+s & 2+8\exp(5)
 \end{bmatrix}'
 \]
 and
 \[
 w=\begin{bmatrix}
 \exp(5)-1 + s & 4+\exp(5)
 \end{bmatrix}'.
 \]
 Calculating  the switching function~$m$
 defined in~\eqref{SwitchingFunction}  yields
 the behavior depicted in Fig.~\ref{fig:m4s}.
 Note that~$m(t)>0$ for $t \in (0,t_1) \cup (t_2,t_3)$,
 and~$m(t)<0$ for~$t \in   (t_1,t_2) \cup (t_3,T)$, so the control~$u$
 satisfies the first-order~MP.

 We now show that the second-order MP
 implies that~$u$ is \emph{not} an optimal control.
 Eq.~\eqref{eq:defhis} yields
\begin{align*}
                        H_1& =A+B,\\
                        H_2& =A- B,\nonumber\\
                        H_3& =\exp(-(A-B)) (A+B) \exp( A-B ),\nonumber\\
                        H_4& =\exp(-(A-B))  \exp(-(A+B))(A-B) \exp(A+B) \exp(A-B).
\end{align*}
Note that
\begin{align}\label{eq:simh3}
                        H_3& =\exp(-(A-B)) (A-B+2B) \exp( A-B )\nonumber\\
                           &=A-B+ 2\exp(-(A-B)) B \exp( A-B ).
\end{align}
%
%
Our goal is to find~$\bar \alpha \in \Q^4$ such that~$r_4(\bar \alpha)>0$. Indeed, this will imply that~$u$ is not optimal.
It turns out that we can find such an~$\bar \alpha$ satisfying~$\bar \alpha_1=1$
 and~$\bar  \alpha_4=0$. Since~$\sum_{i=1}^4 \bar  \alpha_i$ must be zero,
 $\bar  \alpha_3=-1-\bar  \alpha_2$. Then
\begin{align*}
                    \sum_{i=1}^4\bar  \alpha_i H_i  & =   A+B +\bar  \alpha_2(A-B)-( 1+\bar \alpha_2) (A-B+ 2\exp(-(A-B)) B \exp( A-B ))\\
                                               & =2   B  -2(1+\bar  \alpha_2)     \exp(-(A-B)) B \exp( A-B  ),
\end{align*}
so
\begin{align*}
                    \sum_{i=1}^4 \bar  \alpha_i H_i p(t_1)
                                               & =2\left (    B  - (1+\bar \alpha_2)     \exp(-(A-B)) B \exp( A-B  )\right )\exp(A+B)v\\
                                                     & =2 B \exp(A+B)v   -2\sqrt{\rho} (1+\bar \alpha_2)     \exp(-(A-B)) B
                                                      v.
\end{align*}
A tedious but straightforward calculation shows that
\[
                B\exp(A+B)v=\exp(-5)\exp(-(A-B))Bv,
\]
so~$\sum_{i=1}^4 \bar  \alpha_i H_i p(t_1) =0$
for
\[
\bar \alpha_2= (   \sqrt{\rho} \exp(5)   )^{-1} -1 .
\]
Summarizing,~$\bar \alpha=\begin{bmatrix}\bar \alpha_0&  1& (   \sqrt{\rho} \exp(5)   )^{-1}-1 & -(   \sqrt{\rho} \exp(5)   )^{-1}  & 0 \end{bmatrix}' \in \Q^4$.
The second-order term is
\begin{align*}
r_4(\bar \alpha) &= q'(t_1) \sum_{1\leq i < j \leq 4 } \bar \alpha_i \bar \alpha_j [H_i,H_j] p(t_1)\\
                 &= q'(t_1) \left(      \bar \alpha_2 [H_1,H_2]
                 +   \bar \alpha_3 [H_1,H_3]
                 +\bar \alpha_2 \bar \alpha_3 [H_2,H_3]
                 \right)   p(t_1) ,
\end{align*}
and a calculation yields
\[
r_4(\bar \alpha)=\frac{1050 \exp(5)  \left( \exp(5) -1\right) \left(12
   \left(s -3\right)+\exp(5)
   \left(-67+s+\exp(5)\left(67+36
   \exp(5)+12 s\right)\right)\right)}{\left(4+\exp(5) \right)
   \left(1+4 \exp(5) \right) \left(9-3 s+\exp(5) \left(7+9 \exp(5)
   +3 s\right)\right)^2}.
\]
Clearly,~$r_4(\bar \alpha)>0$, so the second-order MP implies that~$u$ in~\eqref{eq:opcande}
is not an optimal control. The reason that~$u$ here 
satisfies the conditions of the first-order MP is that it actually \emph{minimizes}
the spectral radius at time~$T=4$.  
Thus, the second-order MP plays here a similar role to the second-derivative of a function:
it allows to distinguish  between a maximum point and a minimum point.~\qed
 \end{Example}

\section{Proofs}\label{sec:proofs}

\subsection{Proof of Thm.~\ref{thm:secsing}}
Assume that~$u^*(t)\equiv 0$ is an optimal control.
The corresponding solution of~\eqref{eq:transzero} is~$C^*(T)=\exp(AT)$.
For~$\epsilon>0$, consider the control
\be\label{eq:tildeueps}
\tilde{u}(t):= \begin{cases}
            0 ,  &                                    t \in [0, T-4\epsilon^{1/3}), \\
            -1 ,  &                                    t \in [T-4\epsilon^{1/3} ,T- 3 \epsilon^{1/3}), \\
             1 ,  &                                    t \in [T- 3 \epsilon^{1/3}  ,T-   \epsilon^{1/3}), \\
            -1 ,  &                                    t \in [T-   \epsilon^{1/3}  , T). \\
\end{cases}
\ee
Then
\begin{align*}
            C(T,\tilde{u})&= \exp((A-B)\epsilon^{1/3})\exp((A+B)2 \epsilon^{1/3}) \exp((A-B)\epsilon^{1/3})\exp(A (T-4\epsilon^{1/3}))\\
                          &= \exp((A-B)\epsilon^{1/3})\exp((A+B)2 \epsilon^{1/3})  \exp((A-B)\epsilon^{1/3})\exp(   -4\epsilon^{1/3} A) C^*(T),
\end{align*}
and it follows from the computation in~\cite[p. 719]{hermes78} (see also~\cite{krener})  that
\be\label{eq:e13}
            C(T,\tilde{u})= \exp(\frac{2\epsilon}{3}     [B,[B,A]]  )C^*(T)+o(\epsilon).
\ee
Note that this implies  that any result
derived using~$\tilde{u}$ will be  a high-order MP, as
the width of the needle variations in~\eqref{eq:tildeueps} is of order~$\epsilon^{1/3}$
yet the perturbation in~$C(T,\tilde{u})$ with respect to~$C^*(T)$
is of order~$\epsilon$.
By~\eqref{eq:e13},
\[
                        \frac{d}{d \epsilon} C(T,\tilde{u}) |_{\epsilon=0} =  (2/3) [B,[B,A]]  C^*(T)  ,
\]
so
\begin{align*}
            \rho(C(T,\tilde{u})&)-\rho(C^*(T))\\& = (2\epsilon/3) (w^*)' [B,[B,A]]    C^*(T) v^* +o(\epsilon)\\
                                          & = (2\epsilon/3) (w^*)' [B,[B,A]]   \rho(C^*(T))v^* +o(\epsilon).
\end{align*}
If~$(w^*)' [B,[B,A]]      v^* > 0$ then~$\rho(C(T,\tilde{u})) > \rho(C^*(T))$  for all sufficiently
 small~$\epsilon>0$, and this contradicts the optimality of~$u^*$.
This proves~\eqref{eq:abb}.~\qed

\subsection{Proof of Theorem~\ref{thm:mainhere}}
The proof is based on introducing a new control defined by
a
perturbation of the switching times~$ t_0,\dots,t_k$ to
$\tilde{t}_0 :=t_0+s \alpha_0$, $\tilde{t}_1: =t_1+s
(\alpha_0+\alpha_1)$,
 \dots,
  $\tilde{t}_{k }: =t_{k }+s
(\alpha_0+\alpha_1+ \dots+ \alpha_{k })$.
 Here, $s\in \R$ and~$ \alpha \in
\P^{k }$.
Define~$ \tilde{u}(t;s, \alpha)$ by
$\tilde{u}(t)=r$ for $t\in (\tilde{t}_0,\tilde{t}_1)$, $\tilde{u}
(t)=-r$ for $t\in ( \tilde{t}_1,\tilde{t}_2)$, and  so on.
Note that~\eqref{eq:defpk} implies that
the time length of the perturbed control is
\[
\tilde{t}_k-\tilde{t}_0= t_{k }+s
(\alpha_1+ \dots \alpha_{k })-t_0 = t_k-t_0 =T.
\]
 Denote the  corresponding transition matrix
by~$\tilde{ C}(t)=\tilde{ C}(t;s, \alpha)$.
Note also that~$\tilde{u}(\cdot;0,\alpha)=u^*(\cdot)$ for any~$\alpha$, so~$\tilde{C}(\cdot;0,\alpha)=C^*(\cdot)$.
Our goal is to derive an expression
 for  the difference~$ e(s,\alpha):=\rho(\tilde{C}(T;s,\alpha))-\rho(C^*(T))$
  in the form
  \be\label{eq:deffsecorder}
          e(s,\alpha) =s  z_1(\alpha)  +\frac{s^2}{2} z_2(\alpha) +o(s^2),
  \ee
where~$o(s^2) $ denotes a function~$f $ that
satisfies~$\lim_{s\to 0} \frac{f(s)}{s^2}=0$.

Suppose for a moment that~$z_1(\alpha)>0$ [$z_1(\alpha)<0$] for some~$\alpha \in \P^{k }$.
 Then for any sufficiently small~$s>0$ [$s<0$],~\eqref{eq:deffsecorder} implies that
$\rho(C (T;s,\alpha) )  > \rho(C^*(T))$. This contradicts the optimality of~$u^*$, so
\be \label{eq:pisz}
z_1(\alpha)=0,\quad \text{for all } \alpha \in \P^{k }.
\ee
 Thus,~$  e(s,\alpha)=\frac{s^2}{2} z_2(\alpha) +o(s^2),$
and a similar argument implies that the second-order term must satisfy~$z_2(\alpha)\leq 0$.
As we will see below, these   conditions lead to the first- and second-order optimality conditions~\eqref{eq:firsth1} and~\eqref{eq:secoopt}.

The calculation of the terms~$z_1$ and~$z_2$ in~\eqref{eq:deffsecorder} requires two steps.
The first   is to derive an expression for the first- and second-order derivative of~$\tilde{C}(T;s,\alpha)$ with respect to~$s$.
This is based on the Agrachev-Gamkrelidze second-order variation for bang-bang controls~\cite{agrachev-sigalotti} (see also~\cite{ratmansky,yoav}).
 The second step
is to derive an  expression for the  first- and second-order
derivatives of the spectral radius of a matrix   with respect to perturbations of the matrix entries.
This follows the approach  in~\cite{Deutsch19841}.

\subsection{First- and second-order derivatives of the transition matrix}
From here on we
consider
the case where~$k$ is even. (The derivations in the case where~$k$ is odd are similar.)
Thus,
\be \label{eq:xt}
 \tilde{ C}(T;s,\alpha)=
                                \exp( \tilde{\tau}_{k}  Q) \exp (\tilde{\tau}_{k-1}  P)
\dots
             \exp (\tilde{\tau}_2  Q ) \exp (\tilde{\tau}_1  P)  ,
\ee
where~$\tilde{\tau}_i: =\tilde{t}_i-\tilde{t}_{i-1}=\tau_i+s
\alpha_i$. Let~$\dot{  \tilde{C}}: =\frac{d}{ds} \tilde{C} $.

\begin{Proposition}\label{prop:ijuas}
The first- and second-order derivatives of~$  \tilde{C}(T;s,\alpha )$ with respect to~$s$ satisfy
\begin{align}  \label{eq:ghis}
    \tilde{H}  \dot{  \tilde{C}}    &=   \sum_{i=1}^{k }\alpha_i \tilde{G}_i, \nonumber \\
 \tilde{H}  \ddot{ \tilde{ C}}&=
    -  \dot{\tilde{H}}  \dot{ \tilde{ C}}   + \sum_{1\leq i< j \leq k } \alpha_i \alpha_j [ \tilde{G}_i, \tilde{G}_j],
\end{align}
where
\begin{align}\label{eq:defgis}
\tilde{G}_1&:=  \exp(-\tilde{\tau}_1 P)P\exp(\tilde{\tau}_1 P)=P,\nonumber\\
 \tilde{G}_2&:= \exp(-\tilde{\tau}_1 P) Q  \exp(\tilde{\tau}_1 P),\nonumber\\
\tilde{G}_3&:= \exp(-\tilde{\tau}_1 P) \exp(-\tilde{\tau}_2 Q) P  \exp( \tilde{\tau}_2 Q) \exp(\tilde{\tau}_1 P) ,\\
    &\vdots\nonumber
 \end{align}
and
$\tilde{H} = \tilde{H}(s, \alpha) := \exp(-\tilde{\tau}_1  P)      \exp(-\tilde{\tau}_2  Q)
\dots    \exp(-\tilde{\tau}_{k }  Q)   $.
\end{Proposition}

{\sl Proof.}
Rewrite~\eqref{eq:xt} as
\[
 I =     \exp(-\tilde{\tau}_1  P)    \exp(-\tilde{\tau}_2  Q)  \dots
  \exp(-\tilde{\tau}_{k-1}  P)
    \exp(-\tilde{\tau}_{k }  Q)
     \tilde{ C} .
\]
 Differentiating  both sides with respect to~$s$
  and rearranging terms  yields
\begin{align} \label{eq:dds}
\tilde{H} \dot{ \tilde{C}}  & =
     \alpha_1  P  \tilde{C}_0  \\
     & +\alpha_2  \exp(-\tilde{\tau}_1  P)
		     Q \tilde{ C}_1
\nonumber \\
    & +\alpha_3 \exp(-\tilde{\tau}_1  P) \exp(-\tilde{\tau}_2  Q)
      P  \tilde{ C}_2
 \nonumber \\
& \;  \vdots \nonumber \\
&  + \alpha_{k } \exp(-\tilde{\tau}_1  P) \exp(-\tilde{\tau}_2  Q )
\dots \exp(-\tilde{\tau}_{k-1}  P)
    Q   \tilde{ C}_{k-1}  ,   \nonumber
\end{align}
where~$\tilde{ C}_0 : =I$, $\tilde{ C}_1 : =\tilde{C}( \tilde{\tau}_1)$,   $\tilde{ C}_2=\tilde{C}(\tilde{\tau}_1+\tilde{\tau}_2)$, and so on.
This yields the first equation in~\eqref{eq:ghis}.
Differentiating~\eqref{eq:dds} with respect to~$s$
   yields
\begin{align*}
\dot{\tilde{H}}    \dot{ \tilde{ C}}  + \tilde{H}
\ddot{ \tilde{ C}} =
  \sum_{i=1}^{k-1} \sum_{j=i+1} ^k \alpha_i  \alpha_{j } [   \tilde{G}_i,  \tilde{G}_{j }],
\end{align*}
 and this completes the proof of Prop.~\ref{prop:ijuas}.
\qed

 The next step is to determine the first- and second-order derivatives of~$
 \tilde{\rho}(s,\alpha)= \rho(\tilde{C} (T;s,\alpha))$ with respect to~$s$.

\subsection{First- and second-order derivatives  of the spectral radius}

  Let~$\tilde{v}=\tilde{v}(s,\alpha)$, $\tilde{w}=\tilde{w}(s,\alpha)$ denote a nonnegative right and  a
  left eigenvector of~$ \tilde{C} (T)$
 corresponding to the eigenvalue~$\tilde{\rho}$, and
 normalized so that~$ \tilde{w}'\tilde{v}=1$. Note that
 since~$\rho^*$ is simple, 
 the spectral radius of~$\tilde C(T)$ 
 will also be simple for all~$|s|$ sufficiently small. 
  For a matrix~$D$,
 let~$ {D}^ \#$ denote the   Drazin  inverse of~$D$.

 \begin{Proposition}\label{prop:ehyn}
 The first- and second-order derivatives of~$\tilde{\rho}$ with respect to~$s$ satisfy
 \begin{align}\label{eq:rhofp}
                \dot{\tilde{\rho}} |_{s=0}&=     (w^*)'\left( \dot{\tilde{C}} |_{s=0}\right) v^*,\\
                        \ddot{\tilde{\rho}}|_{s=0} & =  (w^*)' \left( \ddot{\tilde{C}} |_{s=0} \right )
                        v^* \nonumber\\&+ 2  (w^*)' \left(\dot{\tilde{C}} |_{s=0} \right)   (\rho^*I - C^*)^ \#  \left(\dot{\tilde{C}}|_{s=0} \right) v^*  .\nonumber
 \end{align}
 \end{Proposition}
{\sl Proof.} Differentiating   the equation
$
                             \tilde{C} \tilde{v}=\tilde{\rho}  \tilde{v}
 $
  with respect to~$s$ yields
 \be\label{eq:firstdertemp}
\dot{ \tilde{C}}{\tilde{v}} +                             \tilde{C}  \dot{\tilde{v} }  =\dot{\tilde{\rho}}  \tilde{v}  + \tilde{\rho}    \dot{\tilde{v}}  .
 \ee
 Multiplying on the left by~$\tilde{w}'$ and using the fact that~$\tilde{w}' \tilde{C}=
 \tilde{\rho}\tilde{w}'$ yields
  \[
                          \tilde{ w}' \dot{\tilde{C}} \tilde{v} +\tilde{ \rho}
                           \tilde{w}' \dot{\tilde{v} }=\dot{\tilde{\rho}} + \tilde{\rho} \tilde{  w}' \dot{\tilde{v}} ,
 \]
 so
 \be\label{eq:der1}
                \dot{\tilde{\rho}}=    \tilde{ w}' \dot{\tilde{C}} \tilde{v}.
 \ee
This proves~\eqref{eq:rhofp}.
To calculate the second-order derivative,   differentiate~\eqref{eq:firstdertemp} with respect to~$s$.
 This yields
 \[
                          \ddot{\tilde{C}} \tilde{v} + 2 \dot{\tilde{C}} \dot{\tilde{v} }  +\tilde{ C}\ddot{\tilde{v }} = \ddot{\tilde{\rho}} \tilde{ v}  +2\dot{\tilde{\rho}}  \dot{\tilde{v}}+\tilde{ \rho}    \ddot{\tilde{v}} .
 \]
Multiplying on the left by~$\tilde{w}'$   yields
  \[
                        \tilde{w}'  \ddot{\tilde{C}} \tilde{v }+ 2 \tilde{w}' \dot{\tilde{C}} \dot{\tilde{v} }  +\tilde{\rho}\tilde{ w}'   \ddot{\tilde{v} } = \ddot{\tilde{\rho}}    +2\dot{\tilde{\rho}} \tilde{w}' \dot{\tilde{v}}+
                        \tilde{\rho} \tilde{ w}'   \ddot{\tilde{v}} ,
 \]
 so
 \begin{align} \label{eq:tempder2}
 \ddot{\tilde{\rho}} &   =                   \tilde{  w}'  \ddot{\tilde{C}}\tilde{ v}
  + 2 \tilde{w}' \dot{\tilde{C}} \dot{\tilde{v} }  -2\dot{\tilde{\rho}}\tilde{w}' \dot{\tilde{v}} \nonumber\\
             &   =                     \tilde{w}'  \ddot{\tilde{C}}\tilde{ v}
              + 2 \tilde{w}' \dot{\tilde{C}} \dot{\tilde{v} }  -2\left( \tilde{w}' \dot{\tilde{C}}
              \tilde{v} \right) \tilde{w}' \dot{\tilde{v}}\nonumber \\
              &   =                     \tilde{w}'  \ddot{\tilde{C}}\tilde{ v}
              + 2 \tilde{w}' \dot{\tilde{C}}  \left(   I-\tilde{v}\tilde{w}'     \right ) \dot{\tilde{v}},
 \end{align}
 where the second equation follows from~\eqref{eq:der1}.
 To simplify this expression, let~$\tilde{D}:=\tilde{\rho} I-\tilde{C}$. It follows from~\eqref{eq:firstdertemp} that
 \be\label{eq:somlp}
   \tilde{D}   \dot{\tilde{v}}=                      \dot{\tilde{C}}  \tilde{v}  - \dot{\tilde{\rho}}  \tilde{v}  .
 \ee
 Since~$\rho^*$ is a simple eigenvalue of~$C^*$,
  $\tilde{\rho}$ is a simple eigenvalue of~$\tilde{C}$ for all~$|s|$ sufficiently small, so
   zero is a \emph{simple} eigenvalue of~$\tilde{D}$. Hence, there exists a unique  generalized inverse~$\tilde{D}^ \#$ of~$\tilde{D}$ satisfying
 \be \label{eq:geninv}
                           \tilde{D}\tilde{D}^ \# =\tilde{D}^ \# \tilde{D}             , \;\;
                           \tilde{D} \tilde{D}^ \# \tilde{D}= \tilde{D}, \;\;
                           \tilde{D}^ \# \tilde{D} \tilde{D}^ \#=\tilde{D}^ \#.
 \ee
 Multiplying~\eqref{eq:somlp} on the left by~$\tilde{D}^ \# $
 and using the fact that~$\tilde{D}^ \# \tilde{v} =0$ and~\eqref{eq:geninv} yields
 \[
                 \tilde{D}  \tilde{D}^ \#  \dot{\tilde{v}}=
                                   \tilde{D}^ \#  \dot{\tilde{C}} \tilde{ v}    .
 \]
 It is not difficult to show that~$I- \tilde{D}  \tilde{D}^ \# =\tilde{v} \tilde{w}'$,
  so
 \[
                (I-\tilde{v}\tilde{w}') \dot{\tilde{v}}=                    \tilde{D}^ \#  \dot{\tilde{C}}  \tilde{v}    .
 \]
 Multiplying this on the left by~$\tilde{w}'\dot{\tilde{C}}$ yields
  \[
             \tilde{w}'\dot{\tilde{C}}  (I-\tilde{v}\tilde{w}') \dot{\tilde{v}}=
                           \tilde{w}'\dot{\tilde{C}}     \tilde{D}^ \#  \dot{\tilde{C}} \tilde{ v}   ,
 \]
 and substituting this in~\eqref{eq:tempder2} yields~$\ddot{\tilde{\rho}}  = \tilde{ w}'  \ddot{\tilde{C}} \tilde{v} + 2 \tilde{ w}'\dot{\tilde{C}}     \tilde{D}^ \#  \dot{\tilde{C}}\tilde{  v}$.
 Setting~$s =0$ completes the proof of Prop.~\ref{prop:ehyn}.~\qed

We can now prove Thm.~\ref{thm:mainhere}. Combining
 \eqref{eq:deffsecorder},  \eqref{eq:ghis}, and~\eqref{eq:rhofp} yields
\[
                z_1(\alpha)= \sum_{i=1}^{k } \alpha_i (w^*)'(H(0,\alpha) )^{-1}G_i v^*,
\]
where~$G_i:=\tilde{G}_i |_{s=0}$. It follows from~\eqref{eq:defgis} and~\eqref{eq:firsth1}  that~$G_i= \exp(-\tau_1 P)H_i\exp(\tau_1 P)$.
The definition of~$q$ in~\eqref{Qeq} implies that
\be\label{eq:wth}
        (w^*)'= q'(T)=q'(0)H(0,\alpha),
\ee
so
\begin{align*}
            z_1(\alpha)&=  \sum_{i=1}^{k }  \alpha_i  q'(0) G_i p(0)\\
                    &= \sum_{i=1}^{k } \alpha_i  q'(t_1) H_i p(t_1).
\end{align*}
Combining this with~\eqref{eq:pisz} proves~\eqref{eq:firsth1}. Note that the proof so far used only first-order derivatives with respect to~$s$.

To prove~\eqref{eq:secoopt}, fix an arbitrary~$\alpha \in \Q^{k }$. Then
by~\eqref{eq:defpkpk},~$ \sum_{i=1}^{k } \alpha_i H_i p(t_1) =0$
and this implies that~$  \sum_{i=1}^{k } \alpha_i G_i p(0)=0  $, so   Prop.~\ref{prop:ijuas} yields
$
           \tilde{ H}(0,\alpha)\left( \dot{\tilde {C}} |_{s=0}\right) v^* = 0,
$
i.e.,
\be   \label{eq:notosub}
                          \left( \dot{\tilde {C}} |_{s=0}\right) v^* = 0.
\ee
Substituting this in~\eqref{eq:ghis} yields
\begin{align*}
 \tilde{H} (0,\alpha)\left( \ddot{ \tilde{ C}}|_{s=0}\right)  v^*&=
      \sum_{1\leq i< j \leq k } \alpha_i \alpha_j [ {G}_i,  G _j]v^*,
\end{align*}
 and multiplying on the left by~$q'(0)$ gives
 \[
         (w^*)'\left( \ddot{ \tilde{ C}}|_{s=0}\right)  v^* =
      q'(0) \sum_{1\leq i< j \leq k } \alpha_i \alpha_j [ {G}_i,  G _j]v^*.
 \]
On the other-hand, substituting~\eqref{eq:notosub}   in~\eqref{eq:rhofp} yields
\begin{align*}
%
                        \ddot{\tilde{\rho}} |_{s=0}  & =  (w^*)' \left( \ddot{\tilde{C}} |_{s=0} \right )
                        v^*  ,
\end{align*}
so
\begin{align*}
            z_2(\alpha)& = q'(0) \sum_{1\leq i< j \leq k } \alpha_i \alpha_j [ {G}_i,  G _j]p(0)\\
            & = q'(t_1) \sum_{1\leq i< j \leq k } \alpha_i \alpha_j [ H_i,  H _j]p(t_1)
\end{align*}
and this proves~\eqref{eq:secoopt}.
This completes the proof of Thm.~\ref{thm:mainhere}.~\qed

\section{Conclusions}
The GUAS problem for switched systems is difficult because of the huge number
of possible switching laws. This led to the variational approach that
is based on trying to characterize   a ``most destabilizing'' switching law.

For a PLSS,  every possible switching law generates
a positive transition matrix and the problem of finding the 
``most destabilizing'' switching law is equivalent to 
finding the switching law that maximizes the spectral radius
of the transition matrix. 

In the  relaxed version of a PLSS, i.e.
a PBCS, this yields a  well-defined optimal control problem, namely,  
for a fixed~$T>0$, find a control~$u$
that maximizes the spectral radius
of the transition matrix  at time~$T$. A first-order necessary condition
for optimality
has been derived in~\cite{lior_SIAM}.
In this paper, we derived a high-order
necessary optimality condition
for a singular control  and for a  bang-bang control. We demonstrated these
conditions using simple examples. 
We are currently trying to use these  high-order conditions
to derive new
theoretical results on the structure of the optimal control in specific problems.
The main technical difficulty is that in general  the Perron right and left
eigenvectors~$v^*$ and~$w^*$ that appear in the optimality conditions 
are complicated functions of the entries of the transition matrix.

\section{Appendix}
\emph{Proof of Thm.~\ref{thm:starho}.}
Pick~$t\geq 0$. Let $\Sigma^t:=\{C(t,u) : u \in \mathcal{U}\}$, i.e., the set of all possible transition matrices at time~$t$.
The definition of~$\U$ implies that~$\Sigma^t$ is  compact.
Note that~\eqref{eq:rhot} yields
\[
\rho(A,B)=\limsup_{t\to \infty }\rho_t(A,B)=\limsup_{t\to \infty } \max_{M \in \Sigma^t} (\rho(M))^{1/t}.
\]
 Pick a multiplicative matrix norm~$||\cdot||:\R^{n\times n} \to \R_+$. 
By the 
\emph{Joint Spectral Radius Theorem} (see, e.g.,~\cite[Ch.~2]{Jungers}),
 \begin{align*} 
 \rho(A,B)
 &=\limsup_{t \rightarrow \infty} \{ ||M||^{1/t} : M \in \Sigma^t \}\\
 &=\limsup_{t \rightarrow \infty} \{ ||C(t,u)||^{1/t} : u \in \U \}.
\end{align*}
If~$\rho (A,B) <1$ 
then
\[
 \limsup_{t \rightarrow \infty} \{ ||C(t,u)||  : u \in \U \} =0,
\]
so~$\lim_{t \rightarrow 0} C(t,u) = 0$ for all~$u \in \U$,
and this  implies GAS of the PBCS~\eqref{eq:pscon}.

Now suppose that~$\rho(A,B) \geqslant 1$.  
By~\cite[p.~22, Corollary~1.2]{Jungers},
there exists a sequence~$R_1,R_2,\dots$ in~$\Sigma^t$ such that 
\[
                        \lim_{k\to\infty} ||  R_k\dots R_1 ||^{1/k} \geq 1. 
\]
In other words, there exists  a sequence of controls~$u_i  \in \U$
such that
\[
                        \lim_{k\to\infty} ||  C(t,u_k) \dots C(t,u_1) ||^{1/k} \geq 1.
\]
Note that~$C(t,u_k) \dots C(t,u_1)$ is the transition matrix at time~$kt$
for the  control obtained by concatenating all the~$u_i$s. 
Thus,  the PBCS is not GAS. This completes the proof.~\QED

 \bibliographystyle{IEEEtran}
 \bibliography{gal_sec_order}

\begin{thebibliography}{10}
\providecommand{\url}[1]{#1}
\csname url@samestyle\endcsname
\providecommand{\newblock}{\relax}
\providecommand{\bibinfo}[2]{#2}
\providecommand{\BIBentrySTDinterwordspacing}{\spaceskip=0pt\relax}
\providecommand{\BIBentryALTinterwordstretchfactor}{4}
\providecommand{\BIBentryALTinterwordspacing}{\spaceskip=\fontdimen2\font plus
\BIBentryALTinterwordstretchfactor\fontdimen3\font minus
  \fontdimen4\font\relax}
\providecommand{\BIBforeignlanguage}[2]{{%
\expandafter\ifx\csname l@#1\endcsname\relax
\typeout{** WARNING: IEEEtran.bst: No hyphenation pattern has been}%
\typeout{** loaded for the language `#1'. Using the pattern for}%
\typeout{** the default language instead.}%
\else
\language=\csname l@#1\endcsname
\fi
#2}}
\providecommand{\BIBdecl}{\relax}
\BIBdecl

\bibitem{sec_order_cdc}
G.~Hochma and M.~Margaliot, ``Stability analysis of positive bilinear control
  systems: A variational approach,'' in \emph{IEEE 52nd Conference on Decision
  and Control (CDC 2013)}, Florence, Italy, 2013, pp. 1355--1359.

\bibitem{lior_SIAM}
L.~Fainshil and M.~Margaliot, ``A maximum principle for positive bilinear
  control systems with applications to positive linear switched systems,''
  \emph{SIAM J.\ Control Optim.}, vol.~50, pp. 2193--2215, 2012.

\bibitem{angeli-sontag-positive-2009}
D.~Angeli, P.~\mbox{de} Leenheer, and E.~D. Sontag, ``Chemical networks with
  inflows and outflows: A positive linear differential inclusions approach,''
  \emph{Biotechnology Progress}, vol.~25, pp. 632--642, 2009.

\bibitem{shorten}
R.~Shorten, F.~Wirth, O.~Mason, K.~Wulff, and C.~King, ``Stability criteria for
  switched and hybrid systems,'' \emph{SIAM Review}, vol.~49, pp. 545--592,
  2007.

\bibitem{branicky_98}
M.~S. Branicky, ``Multiple {L}yapunov functions and other analysis tools for
  switched and hybrid systems,'' \emph{IEEE Trans.\ Automat.\ Control},
  vol.~43, pp. 475--482, 1998.

\bibitem{decarloetal}
R.~DeCarlo, M.~Branicky, S.~Pettersson, and B.~Lennartson, ``Perspectives and
  results on the stability and stabilizability of hybrid systems,'' \emph{Proc.
  IEEE}, vol.~88, pp. 1069--1082, 2000.

\bibitem{libsur99}
D.~Liberzon and A.~S. Morse, ``Basic problems in stability and design of
  switched systems,'' \emph{IEEE Control Systems Magazine}, vol.~19, pp.
  59--70, 1999.

\bibitem{liberzon_book}
D.~Liberzon, \emph{Switching in Systems and Control}.\hskip 1em plus 0.5em
  minus 0.4em\relax Birkh{\"a}user, 2003.

\bibitem{morse-book}
A.~S. Morse, Ed., \emph{Control Using Logic-Based Switching}.\hskip 1em plus
  0.5em minus 0.4em\relax London: Springer, 1997.

\bibitem{Jungers}
R.~Jungers, \emph{The Joint Spectral Radius: Theory and Applications}, ser.
  Lecture Notes in Control and Information Sciences.\hskip 1em plus 0.5em minus
  0.4em\relax Springer, 2009, vol. 385.

\bibitem{johansson}
M.~Johansson, \emph{Piecewise Linear Control Systems}, ser. Lecture Notes in
  Control and Information Sciences.\hskip 1em plus 0.5em minus 0.4em\relax
  Springer-Verlag, 2003, vol. 284.

\bibitem{sun_ge_book}
Z.~Sun and S.~S. Ge, \emph{Switched Linear Systems: Control and Design}.\hskip
  1em plus 0.5em minus 0.4em\relax Springer, 2005.

\bibitem{sun_ge_2}
Z.~Sun and S.~S. Ge, \emph{Stability Theory of Switched Dynamical
  Systems}.\hskip 1em plus 0.5em minus 0.4em\relax Springer, 2011.

\bibitem{farina2000}
L.~Farina and S.~Rinaldi, \emph{Positive Linear Systems: Theory and
  Applications}.\hskip 1em plus 0.5em minus 0.4em\relax John Wiley, 2000.

\bibitem{berman87}
A.~Berman and R.~J. Plemmons, \emph{Nonnegative Matrices in the Mathematical
  Sciences}.\hskip 1em plus 0.5em minus 0.4em\relax SIAM, 1987.

\bibitem{posi_sys89}
M.~A. Krasnoselskij, J.~A. Lifshits, and A.~V. Sobolev, \emph{Positive Linear
  Systems: the Method of Positive Operators}.\hskip 1em plus 0.5em minus
  0.4em\relax Heldermann-Verlag, 1989.

\bibitem{mason-shorten03}
O.~Mason and R.~Shorten, ``A conjecture on the existence of common quadratic
  {L}yapunov functions for positive linear systems,'' in \emph{Proc.\ 2003
  American Control Conf.}, Denver, CO, 2003, pp. 4469--4470.

\bibitem{Gurvits_Olshevsky}
L.~Gurvits and A.~Olshevsky, ``On the {NP}-hardness of checking matrix polytope
  stability and continuous-time switching stability,'' \emph{IEEE Trans.\
  Automat.\ Control}, vol.~54, no.~2, pp. 337--341, 2009.

\bibitem{gurvits-shorten-mason07}
L.~Gurvits, R.~Shorten, and O.~Mason, ``On the stability of switched positive
  linear systems,'' \emph{IEEE Trans.\ Automat.\ Control}, vol.~52, pp.
  1099--1103, 2007.

\bibitem{gurvits-cdc03}
L.~Gurvits, ``What is the finiteness conjecture for linear continuous time
  inclusions?'' in \emph{Proc.\ 42nd IEEE Conf. on Decision and Control}, Maui,
  HI, 2003, pp. 1165--1169.

\bibitem{copos}
O.~Mason, V.~S. Bokharaie, and R.~Shorten, ``Stability and {D}--stability for
  switched positive systems,'' in \emph{Positive Systems}, ser. Lecture Notes
  in Control and Information Sciences, R.~Bru and S.~Romero-Vivo, Eds.\hskip
  1em plus 0.5em minus 0.4em\relax Springer-Verlag, 2009, vol. 389, pp.
  101--109.

\bibitem{forna}
E.~Fornasini and M.~E. Valcher, ``Linear copositive \mbox{L}yapunov functions
  for continuous-time positive switched systems,'' \emph{IEEE Trans.\ Automat.\
  Control}, vol.~55, pp. 1933--1937, 2010.

\bibitem{ron_margaliot}
\BIBentryALTinterwordspacing
O.~Ron, M.~Margaliot, and M.~S. Branicky, ``Switching between linear consensus
  protocols: a variational approach,'' submitted. [Online]. Available:
  \url{http://arxiv.org/abs/1407.2399}
\BIBentrySTDinterwordspacing

\bibitem{mar-bra-full}
M.~Margaliot and M.~S. Branicky, ``Nice reachability for planar bilinear
  control systems with applications to planar linear switched systems,''
  \emph{IEEE Trans.\ Automat.\ Control}, vol.~54, pp. 1430--1435, 2009.

\bibitem{lior}
L.~Fainshil, M.~Margaliot, and P.~Chigansky, ``On the stability of positive
  linear switched systems under arbitrary switching laws,'' \emph{IEEE Trans.\
  Automat.\ Control}, vol.~54, pp. 897--899, 2009.

\bibitem{pyat70}
E.~S. Pyatnitskii, ``Absolute stability of nonstationary nonlinear systems,''
  \emph{Automat. Remote Control}, vol.~1, pp. 5--15, 1970.

\bibitem{pyat71}
E.~S. Pyatnitskii, ``Criterion for the absolute stability of second-order
  nonlinear controlled systems with one nonlinear nonstationary element,''
  \emph{Automat. Remote Control}, vol.~1, pp. 5--16, 1971.

\bibitem{mar-simple}
M.~Margaliot, ``Stability analysis of switched systems using variational
  principles: An introduction,'' \emph{Automatica}, vol.~42, pp. 2059--2077,
  2006.

\bibitem{bar-cdc}
N.~E. Barabanov, ``Lyapunov exponent and joint spectral radius: Some known and
  new results,'' in \emph{Proc.\ 44th IEEE Conf. on Decision and Control},
  Seville, Spain, 2005, pp. 2332--2337.

\bibitem{rap}
L.~B. Rapoport, ``Asymptotic stability and periodic motions of selector-linear
  differential inclusions,'' in \emph{Robust Control via Variable Structure and
  Lyapunov Techniques}, ser. Lecture Notes in Control and Information Sciences,
  F.~Garofalo and L.~Glielmo, Eds.\hskip 1em plus 0.5em minus 0.4em\relax
  Springer, 1996, vol. 217, pp. 269--285.

\bibitem{Boscain2008Nondiagonalizable_case}
M.~Balde and U.~Boscain, ``Stability of planar switched systems: the
  nondiagonalizable case,'' \emph{Communications Pure Applied Analysis},
  vol.~7, pp. 1--21, 2008.

\bibitem{Boscain2009stability_conditions}
M.~Balde, U.~Boscain, and P.~Mason, ``A note on stability conditions for planar
  switched systems,'' \emph{Int. J. Control}, vol.~82, pp. 1882--1888, 2009.

\bibitem{monovich1}
T.~Monovich and M.~Margaliot, ``Analysis of discrete-time linear switched
  systems: A variational approach,'' \emph{SIAM J.\ Control Optim.}, vol.~49,
  pp. 808--829, 2011.

\bibitem{tal2}
T.~Monovich and M.~Margaliot, ``A second-order maximum principle for
  discrete--time bilinear control systems with applications to discrete--time
  linear switched systems,'' \emph{Automatica}, vol.~47, pp. 1489--1495, 2011.

\bibitem{mar-bra-posi}
M.~Margaliot and M.~S. Branicky, ``Stability analysis of positive linear
  switched systems: a variational approach,'' in \emph{Proc. 6th \mbox{IFAC}
  Symposium on Robust Control Design~(ROCOND09)}, Haifa, Israel, 2009.

\bibitem{filippov-paper}
A.~F. Filippov, ``On certain questions in the theory of optimal control,''
  \emph{SIAM J.\ Control Optim.}, vol.~1, pp. 76--84, 1962.

\bibitem{mat_ana_sec_ed}
R.~A. Horn and C.~R. Johnson, \emph{Matrix Analysis}, 2nd~ed.\hskip 1em plus
  0.5em minus 0.4em\relax Cambridge University Press, 2013.

\bibitem{Deutsch19841}
E.~Deutsch and M.~Neumann, ``Derivatives of the {P}erron root at an essentially
  nonnegative matrix and the group inverse of an {M}-matrix,'' \emph{J. Math.
  Anal. Appl.}, vol. 102, no.~1, pp. 1--29, 1984.

\bibitem{agrachev-sigalotti}
A.~A. Agrachev and M.~Sigalotti, ``On the local structure of optimal
  trajectories in {R}$^3$,'' \emph{SIAM J.\ Control Optim.}, vol.~42, pp.
  513--531, 2003.

\bibitem{hermes78}
H.~Hermes, ``Lie algebras of vector fields and local approximation of
  attainable sets,'' \emph{SIAM J.\ Control Optim.}, vol.~16, pp. 715--727,
  1978.

\bibitem{krener}
A.~J. Krener, ``The high order maximal principle and its applications to
  singular extremals,'' \emph{SIAM J.\ Control Optim.}, vol.~15, pp. 256--293,
  1977.

\bibitem{ratmansky}
I.~Ratmansky and M.~Margaliot, ``A simplification of the
  {A}grachev--{G}amkrelidze second-order variation for bang--bang controls,''
  \emph{Systems Control Lett.}, vol.~59, pp. 25--32, 2010.

\bibitem{yoav}
Y.~Sharon and M.~Margaliot, ``Third-order nilpotency, finite switchings and
  asymptotic stability,'' \emph{J. Diff. Eqns.}, vol. 233, pp. 136--150, 2007.

\end{thebibliography}

\end{document}